\documentclass{article}
\usepackage{amsmath, amsthm, amsfonts, times,latexsym,graphicx,color}
\usepackage{url}

\newtheorem{theorem}{Theorem}[section]
\newtheorem{lemma}[theorem]{Lemma}

\newtheorem{corollary}[theorem]{Corollary}
\newtheorem{proposition}[theorem]{Proposition}
\newtheorem{example}[theorem]{Example}
\newtheorem{remark}[theorem]{Remark}

\def\deq{ \ \dot= \ }

\def\Ex{\mathop{\rm I\!E}\nolimits}
\def\Pr{\mathop{\rm I\!P}\nolimits}

\def\KNem{K_{\rm Nem}}
\def\KT2{K_{\rm Type2}}
\def\KTB{K_{\rm TrBern}}

\def\B{\mathbb{B}}
\def\H{\mathbb{H}}

\def\R{\mathbb{R}}

\def\Var{\mathop{\rm Var}\nolimits}

\begin{document}

\title{Nemirovski's Inequalities Revisited}
\author{Lutz D\"umbgen, Sara A. van de Geer, Mark C. Veraar, and Jon A. Wellner}
\date{February 2009}
\maketitle

\section{Introduction.}

Our starting point is the following well known theorem from probability: 
Let $X_1, \ldots , X_n$ be (stochastically) independent random variables with 
finite second moments, and let $S_n = \sum_{i=1}^n X_i$.  Then
\begin{eqnarray}
	\Var(S_n ) = \sum_{i=1}^n \Var(X_i) .
\label{VarianceOfSumEqualsSumOfVariances}
\end{eqnarray}
If we suppose that each $X_i$ has mean zero, $\Ex X_i = 0$, then 
(\ref{VarianceOfSumEqualsSumOfVariances}) becomes
\begin{eqnarray}
	\Ex S_n^2 = \sum_{i=1}^n \Ex X_i^2 .
\label{SecondMomenOfSumEqualsSumOfSecondMoments}
\end{eqnarray}
This equality generalizes easily to vectors in a Hilbert space $\H$ with inner product 
$\langle \cdot , \cdot \rangle$: If the $X_i$'s are independent with values in $\H$ such 
that $\Ex X_i = 0$ and $\Ex \|X_i\|^2 < \infty$, then 
$\|S_n\|^2 = \langle S_n, S_n \rangle = \sum_{i,j=1}^n \langle X_i , X_j \rangle$, and 
since $\Ex \langle X_i , X_j \rangle = 0$ for $ i \not= j$ by independence,
\begin{eqnarray}
	\Ex \| S_n \|^2
	= \sum_{i,j=1}^n \Ex \langle X_i , X_j \rangle
	= \sum_{i=1}^n \Ex \| X_i \|^2 .
\label{HilbertSpaceVersionSecondMomenOfSumEqualsSumOfSecondMoments}
\end{eqnarray}

What happens if the $X_i$'s take values in a (real) Banach space $(\B, \|\cdot\|)$? 
In such cases, in particular when the square of the norm $\| \cdot \|$ is not given by 
an inner product, we are aiming at inequalities of the following type: Let $X_1$, 
$X_2$, \ldots, $X_n$ be independent random vectors with values in $(\B, \|\cdot\|)$ 
with $\Ex X_i = 0$ and $\Ex \|X_i\|^2 < \infty$. With $S_n := \sum_{i=1}^n X_i$ we 
want to show that
\begin{equation}
\label{ineq: Nemirovski}
	\Ex \|S_n\|^2 \ \le \ K \sum_{i=1}^n \Ex \|X_i\|^2
\end{equation}
for some constant $K$ depending only on $(\B,\|\cdot\|)$.

For statistical applications, the case $(\B,\|\cdot\|) = \ell_r^d := (\R^d, \|\cdot\|_r)$ 
for some $r \in [1,\infty]$ is of particular interest. 
Here the $r$-norm of a vector $x \in \R^d$ is defined as
\begin{equation}
\label{eq: Definition of r-norm}
	\|x\|_r \ := \begin{cases}
		\displaystyle
		\Bigl( \sum_{j=1}^d |x_j|^r \Bigr)^{1/r} & \text{if} \ 1 \le r < \infty , \\
		\displaystyle
		\max_{1 \le j \le d} \, |x_j|            & \text{if} \ r = \infty .
	\end{cases}
\end{equation}
An obvious question is how the exponent $r$ and the dimension $d$ enter an inequality 
of type (\ref{ineq: Nemirovski}). The influence of the dimension $d$ is crucial, since 
current statistical research often involves small or moderate ``sample size'' $n$ 
(the number of independent units), say on the order of $10^2$ or $10^4$, while the number 
$d$ of items measured for each independent unit is large, say on the order of $10^6$ or $10^7$. 
The following two examples for the random vectors $X_i$ provide lower bounds for the 
constant $K$ in (\ref{ineq: Nemirovski}):

\begin{example} \textbf{(A lower bound in $\ell_r^d$)} \
\label{ex: lower bound ell_r}
Let $b_1, b_2, \ldots, b_d$ denote the standard basis of $\R^d$, and let 
$\epsilon_1, \epsilon_2, \ldots, \epsilon_d$ be independent \textsl{Rademacher variables}, 
i.e. random variables taking the 
values $+1$ and $-1$ each with probability $1/2$. Define $X_i := \epsilon_i b_i$ 
for $1 \le i \le n := d$. Then $\Ex X_i = 0$, $\|X_i\|_r^2 = 1$, and 
$\|S_n\|_r^2 = d^{2/r} = d^{2/r-1} \sum_{i=1}^n \|X_i\|_r^2$. Thus any candidate for 
$K$ in (\ref{ineq: Nemirovski}) has to satisfy $K \ge d^{2/r - 1}$.
\end{example}

\begin{example} \textbf{(A lower bound in $\ell_\infty^d$)} \
\label{ex: lower bound ell_infty}
Let $X_1, X_2, X_3, \ldots$ be independent random vectors, uniformly 
distributed on $\{-1,1\}^d$ each. Then $\Ex X_i = 0$ and $\|X_i\|_\infty = 1$. 
On the other hand, according to the Central Limit Theorem, $n^{-1/2} S_n$ 
converges in distribution as $n \to \infty$ to a random vector 
$Z = (Z_j)_{j=1}^d$ with independent, standard Gaussian components, $Z_j \sim N(0,1)$. Hence
$$
	\sup_{n \ge 1} \frac{\Ex \|S_n\|_\infty^2}{\sum_{i=1}^n \Ex \|X_i\|_\infty^2}
	= \sup_{n \ge 1} \Ex \|n^{-1/2} S_n\|_\infty^2
	\ \ge \ \Ex \|Z\|_\infty^2
	= \Ex \max_{1 \le j \le d} Z_j^2 .
$$
But it is well-known that $\max_{1 \le j \le d} |Z_j| = \sqrt{2 \log d} + o_p(1)$ 
as $d \to \infty$. Thus candidates $K(d)$ for the constant in (\ref{ineq: Nemirovski}) have to satisfy
$$
	\liminf_{d \to \infty} \frac{K(d)}{2 \log d}
	\ \ge \ 1 .
$$
\end{example}

At least three different methods have been developed to prove inequalities of the 
form given by (\ref{ineq: Nemirovski}). The three approaches known to us are: 
\begin{quote}
(a)  deterministic inequalities for norms;\\
(b)  probabilistic methods for Banach spaces;\\
(c)  empirical process methods.
\end{quote}
Approach (a) was used by Nemirovski \cite{Nemirovski:00} to show that in the 
space $\ell_r^d$ with $d \ge 2$, inequality (\ref{ineq: Nemirovski}) holds with 
$K = C \min(r, \log(d))$ for some universal (but unspecified) constant $C$. 
In view of Example~\ref{ex: lower bound ell_infty}, this constant has the correct 
order of magnitude if $r = \infty$. For statistical applications see Greenshtein and Ritov 
\cite{GreenshteinRitov:04}. Approach (b) uses special moment inequalities from probability 
theory on Banach spaces which involve nonrandom vectors in $\B$ and Rademacher 
variables as introduced in Example~\ref{ex: lower bound ell_r}. Empirical process theory 
(approach (c)) in general deals with sums of independent random elements in 
infinite-dimensional Banach spaces. By means of chaining arguments, metric entropies 
and approximation arguments, ``maximal inequalities'' for such random sums are built 
from basic inequalities for sums of independent random variables or finite-dime
 nsional random vectors, in particular, ``exponential inequalities''; see e.g.\ 
 Dudley \cite{Dudley:99}, van der Vaart and Wellner \cite{vanderVaartWellner:96}, 
 Pollard \cite{Pollard:90}, de la Pena and Gin\'{e} \cite{delaPenaGine:99}, 
 or van de Geer \cite{vandeGeer:00}.

Our main goal in this paper is to compare the inequalities resulting from these 
different approaches and to refine or improve the constants $K$ obtainable by each method. 
The remainder of this paper is organized as follows: In Section~\ref{sec: Geometry} we 
review several deterministic inequalities for norms and, in particular, key arguments of 
Nemirovski \cite{Nemirovski:00}. Our exposition includes explicit and improved constants. 
While finishing the present paper we became aware of yet unpublished work of 
\cite{Nemirovski:04} and \cite{JuditskyNemirovski:08} who also improved some inequalities 
of \cite{Nemirovski:00}. Rio \cite{Rio:08} uses similar methods in a different context. 
In Section~\ref{sec: Type and cotype} we present inequalities of type (\ref{ineq: Nemirovski}) 
which follow from type and co-type inequalities developed in probability theory on Banach spaces. 
In addition, we provide and utilize a new type inequality for the normed space $\ell_\infty^d$. 
To do so we utilize, among other tools, exponential inequalities of Hoeffding \cite{MR0144363} 
and Pinelis \cite{MR1272088}. In Section~\ref{sec: Bernstein} we follow approach (c) and 
treat $\ell_\infty^d$ by means of a truncation argument and Bernstein's exponential inequality.
Finally, in Section~\ref{sec: Comparisons} we compare the inequalities resulting from these three 
approaches. In that section we relax the assumption that $\Ex X_i = 0$ for a more thorough 
understanding of the differences between the three approaches. Most proofs are deferred to 
Section~\ref{sec: Proofs}.

\section{Nemirovski's approach: Deterministic inequalities for\\
	norms.}
\label{sec: Geometry}

In this section we review and refine inequalities of type (\ref{ineq: Nemirovski}) 
based on deterministic inequalities for norms. The considerations for $(\B,\|\cdot\|) = \ell_r^d$ 
follow closely the arguments of \cite{Nemirovski:00}.

\subsection{Some inequalities for $\R^d$ and the norms $\|\cdot\|_r$}

Throughout this subsection let $\B = \R^d$, equipped with one of the norms 
$\|\cdot\|_r$ defined in (\ref{eq: Definition of r-norm}).  For $x \in \R^d$ we think of
$x$ as a column vector and write $x^\top $ for the corresponding row vector.  
Thus $x x^\top$ is a $d \times d$ matrix with entries 
$x_i x_j$ for $i,j \in \{ 1, \ldots , d \}$.

\paragraph{A first solution.} Recall that for any $x \in \R^d$,
\begin{equation}
\label{ineq: norms}
	\|x\|_r \ \le \ \|x\|_q \ \le \ d^{1/q - 1/r} \|x\|_r
	\quad\text{for} \ 1 \le q < r \le \infty .
\end{equation}
Moreover, as mentioned before,
$$
	\Ex \|S_n\|_2^2 \ = \ \sum_{i=1}^n \Ex \|X_i\|_2^2 .
$$
Thus for $1 \le q < 2$,
$$
	\Ex \|S_n\|_q^2
	\ \le \ (d^{1/q - 1/2})_{}^2 \Ex \|S_n\|_2^2
	\ =   \ d^{2/q - 1} \sum_{i=1}^n \Ex \|X_i\|_2^2
	\ \le \ d^{2/q - 1} \sum_{i=1}^n \Ex \|X_i\|_q^2 ,
$$
whereas for $2 < r \le \infty$,
$$
	\Ex \|S_n\|_r^2
	\ \le \ \Ex \|S_n\|_2^2
	\ =   \ \sum_{i=1}^n \Ex \|X_i\|_2^2
	\ \le \ d^{1 - 2/r} \sum_{i=1}^n \Ex \|X_i\|_r^2 .
$$
Thus we may conclude that (\ref{ineq: Nemirovski}) holds with
$$
	K = \widetilde{K}(d,r) \ := \ \begin{cases}
		d^{2/r - 1} & \text{if} \ 1 \le r \le 2 , \\
		d^{1 - 2/r} & \text{if} \ 2 \le r \le \infty .
	\end{cases}
$$
Example~\ref{ex: lower bound ell_r} shows that this constant $\widetilde{K}(d,r)$ is 
indeed optimal for $1 \le r \le 2$.

\paragraph{A refinement for $r > 2$.}
In what follows we shall replace $\widetilde{K}(d,r) = d^{1 - 2/r}$ with substantially 
smaller constants. The main ingredient is the following result:

\begin{lemma}
\label{lem: Nemirovski}
For arbitrary fixed $r \in [2,\infty)$ and $x \in \R^d \setminus \{0\}$ let
$$
	h(x) \ := \ 2 \|x\|_r^{2-r} \bigl( |x_i|^{r-2} x_i \bigr)_{i=1}^d
$$
while $h(0) := 0$. Then for arbitrary $x, y \in \R^d$,
$$
	\|x\|_r^2 + h(x)^\top y
	\ \le \ \|x + y\|_r^2
	\ \le \ \|x\|_r^2 + h(x)^\top y + (r - 1) \|y\|_r^2 .
$$
\end{lemma}

\cite{NemirovskiYudin:83} and \cite{Nemirovski:00} stated Lemma~\ref{lem: Nemirovski} 
with the factor $r-1$ on the right side replaced with $C r$ for some (absolute) constant 
$C > 1$. Lemma~\ref{lem: Nemirovski}, 
which is a special case of the more general Lemma~\ref{lem: good news!} in the next 
subsection, may be applied to the partial sums $S_0 := 0$ and 
$S_k := \sum_{i=1}^k X_i$, $1 \le k \le n$, to show that for $2 \le r < \infty$,
\begin{eqnarray*}
	\Ex \|S_k\|_r^2
	& \le & \Ex \bigl( \|S_{k-1}\|_r^2 + h(S_{k-1})^\top X_k + (r - 1) \|X_k\|_r^2 \bigr) \\
	& = & \Ex \|S_{k-1}\|_r^2 + \Ex h(S_{k-1})^\top \Ex X_k + (r - 1) \Ex \|X_k\|_r^2 \\
	& = & \Ex \|S_{k-1}\|_r^2 + (r - 1) \Ex \|X_k\|_r^2 ,
\end{eqnarray*}
and inductively we obtain a second candidate for $K$ in (\ref{ineq: Nemirovski}):
$$
	\Ex \|S_n\|_r^2 \ \le \ (r - 1) \sum_{i=1}^n \Ex \|X_i\|_r^2
	\quad\text{for} \ 2 \le r < \infty .
$$
Finally, we apply (\ref{ineq: norms}) again: For $2 \le q \le r \le \infty$ with $q < \infty$,
$$
	\Ex \|S_n\|_r^2
	\ \le \ \Ex \|S_n\|_q^2 \\
	\ \le \ (q - 1) \sum_{i=1}^n \Ex \|X_i\|_q^2 \\
	\ \le \ (q - 1) d^{2/q - 2/r} \sum_{i=1}^n \Ex \|X_i\|_r^2 .
$$
This inequality entails our first ($q = 2$) and second ($q = r < \infty$) 
preliminary result, and we arrive at the following refinement:

\begin{theorem}
\label{thm: Nemirovski}
For arbitrary $r \in [2, \infty]$,
$$
	\Ex \|S_n\|_r^2
	\ \le \ \KNem(d,r) \sum_{i=1}^n \Ex \|X_i\|_r^2
$$
with
$$
	\KNem(d,r) \ := \ \inf_{q \in [2,r] \cap \R} (q - 1) d^{2/q - 2/r} .
$$
This constant $\KNem(d,r)$ satisfies the (in)equalities
$$
	\KNem(d,r) \ \left\{\begin{array}{ccl}
		=   & d^{1 - 2/r} & \text{if} \ d \le 7 \\
		\le & r - 1 \\
		\le & 2 e \log d - e & \text{if} \ d \ge 3 ,
	\end{array}\right.
$$
and
$$
	\KNem(d,\infty) \ \ge \ 2 e \log d - 3e .
$$
\end{theorem}

\begin{corollary}
\label{cor: NemirovskiForREqualsInfinity}
In case of $(\B,\|\cdot\|) = \ell_\infty^d$ with $d \ge 3$, inequality~(\ref{ineq: Nemirovski}) 
holds with constant $K = 2 e \log d - e$. If the $X_i$'s are also identically distributed, then
$$
	\Ex \| n^{-1/2} S_n \|_\infty^2
	\ \le \ (2e \log d - e) \Ex \| X_1 \|_\infty^2 .
$$
\end{corollary}

Note that
$$
	\lim_{d \to \infty} \frac{\KNem(d,\infty)}{2 \log d}
	\ = \ \lim_{d \to \infty} \frac{2e \log d - e}{2 \log d}
	\ = \ e .
$$
Thus Example~\ref{ex: lower bound ell_infty} entails that for large dimension $d$, 
the constants $\KNem(d,\infty)$ and $2e \log d - e$ are optimal up to a factor close to $e \deq 2.7183$.

\subsection{Arbitrary $L_r$-spaces}
\label{subsec: L_r-spaces}

Lemma~\ref{lem: Nemirovski} is a special case of a more general inequality: 
Let $(T, \Sigma, \mu)$ be a $\sigma$-finite measure space, and for 
$1 \le r < \infty$ let $L_r(\mu)$ be the set of all measurable functions 
$f : T \to \R$ with finite (semi-) norm
$$
	\|f\|_r \ := \ \Bigl( \int |f|^r \, d\mu \Bigr)^{1/r} ,
$$
where two such functions are viewed as equivalent if they coincide almost 
everywhere with respect to $\mu$. In what follows we investigate the functional
$$
	f \ \mapsto \ V(f) := \|f\|_r^2
$$
on $L_r(\mu)$. Note that $(\R^d, \|\cdot\|_r)$ 
corresponds to $(L_r(\mu), \|\cdot\|_r)$ if we take $T = \{1,2,\ldots,d\}$ equipped 
with counting measure $\mu$.

Note that $V(\cdot)$ is convex;
thus for fixed $f, g \in L_r(\mu)$, the function 
$$
       v(t) \ := \ V(f+tg) = \| f+t g\|_r^2, \qquad t \in \R
$$
is convex with derivative 
$$
       v'(t) = v^{1-r/2} (t) \int 2 | f + t g |^{r-2} (f+tg) g \, d\mu .
$$
By convexity of $v$ it follows that 
$$
V(f+g) - V(f) = v(1) - v(0) \ge v' (0) \ := \ DV (f,g).
$$
This proves the lower bound in the following lemma.  
We will prove the upper bound in Section~\ref{sec: Proofs} 
by computation of $v'' $ and application of H\"older's inequality.

\begin{lemma}
\label{lem: good news!}
Let $r \ge 2$. Then for arbitrary $f, g \in L_r(\mu)$,
$$
	DV(f,g) \ = \ \int h(f) g \, d\mu
	\quad\text{with}\quad
	h(f) \ := \ 2 \| f \|_r^{2-r} |f|^{r-2} f \ \in \ L_q(\mu) ,
$$
where $q := r/(r-1)$. Moreover,
$$
	V(f) + DV(f,g) \ \le \ V(f + g) \ \le \ V(f) + DV(f,g) + (r - 1) V(g) .
$$
\end{lemma}

\begin{remark}
The upper bound for $V(f + g)$ is sharp in the following sense: 
Suppose that $\mu(T) < \infty$, and let $f, g_o : T \to \R$ be measurable 
such that $|f| \equiv |g_o| \equiv 1$ and $\int f g_o \, d\mu = 0$. 
Then our proof of Lemma~\ref{lem: good news!} reveals that
$$
	\frac{V(f + tg_o) - V(f) - DV(f,tg_o)}{V(t g_o)} \ \to \ r - 1
	\quad\text{as} \ t \to 0 .
$$
\end{remark}

\begin{remark}
In case of $r = 2$, Lemma~\ref{lem: good news!} is well known and easily verified. 
Here the upper bound for $V(f + g)$ is even an equality, i.e.
$$
	V(f + g) \ = \ V(f) + DV(f,g) + V(g) .
$$
\end{remark}

\begin{remark}
Lemma~\ref{lem: good news!} improves on an inequality of \cite{NemirovskiYudin:83}. 
After writing this paper we realized Lemma~\ref{lem: good news!} is also proved by 
\cite{MR1331198}; see his (2.2) and Proposition 2.1, page 1680.
\end{remark}

Lemma~\ref{lem: good news!} leads directly to the following result:

\begin{corollary}
\label{cor: Nemirovski 1}
In case of $\B = L_r(\mu)$ for $r\ge 2$, inequality~(\ref{ineq: Nemirovski}) is satisfied with $K = r - 1$.
\end{corollary}

\subsection{A connection to geometrical functional analysis}

For any Banach space $(\B, \|\cdot\|)$ and Hilbert space 
$(\H,\langle\cdot,\cdot\rangle,\|\cdot\|)$, their Banach-Mazur 
distance $D(\B,\H)$ is defined to be the infimum of
$$
	\|T\| \cdot \|T^{-1}\|
$$
over all linear isomorphisms $T : \B \to \H$, where $\|T\|$ and $\|T^{-1}\|$ denote the usual operator norms
\begin{eqnarray*}
	\|T\|      & := & \sup \bigl\{ \|Tx\| : x \in \B, \|x\| \le 1 \bigr\} , \\
	\|T^{-1}\| & := & \sup \bigl\{ \|T^{-1}y\| : y \in \H, \|y\| \le 1 \bigr\} .
\end{eqnarray*}
(If no such bijection exists, one defines $D(\B,\H) := \infty$.) Given such a bijection $T$,
\begin{eqnarray*}
	\Ex \|S_n\|^2
	& \le & \|T^{-1}\|^2 \Ex \| TS_n\|^2 \\
	& = & \|T^{-1}\|^2 \sum_{i=1}^n \Ex \| TX_i\|^2 \\
	& \le & \|T^{-1}\|^2 \|T\|^2 \sum_{i=1}^n \Ex \|X_i\|^2 .
\end{eqnarray*}
This leads to the following observation:

\begin{corollary}
\label{cor: Nemirovski 2}
For any Banach space $(\B, \|\cdot\|)$ and any Hilbert space 
$(\H,\langle,\cdot,\cdot,\rangle,\|\cdot\|)$ with finite Banach-Mazur distance $D(\B,\H)$, 
inequality~(\ref{ineq: Nemirovski}) is satisfied with $K = D(\B,\H)^2$.
\end{corollary}

A famous result from geometrical functional analysis is John's theorem 
(cf.\ \cite{Tomczak-Jaegermann:89}, \cite{JohnsonLindenstrauss:01}) 
for finite-dimensional normed spaces. It entails that 
$D(\B,\ell_2^{\dim(\B)}) \le \sqrt{\dim(\B)}$. This entails the following fact:

\begin{corollary}
\label{cor: Nemirovski 3}
For any normed space $(\B, \|\cdot\|)$ with finite 
dimension, inequality~(\ref{ineq: Nemirovski}) is satisfied with $K = \dim(\B)$.
\end{corollary}

Note that Example~\ref{ex: lower bound ell_r} with $r = 1$ provides an example 
where the constant $K = \dim(\B)$ is optimal.

\section{The probabilistic approach: Type and co-type\\
	inequalities.}
\label{sec: Type and cotype} 

\subsection{Rademacher type and cotype inequalities}

Let $\{\epsilon_i\}$ denote a sequence of independent Rademacher random 
variables. Let $1 \le p < \infty$. A Banach space $\B$ with norm $\|\cdot\|$ is 
said to be of \textsl{(Rademacher) type $p$} if there is a constant $T_p$ such 
that for all finite sequences $\{ x_i \}$ in $\B$,
$$
	\Ex \Bigl\| \sum_{i=1}^n \epsilon_i x_i \Bigr\|^p
	\ \le \ T_p^p \sum_{i=1}^n \|x_i\|^p .
$$
Similarly, for $1 \le q < \infty$, $\B$ is of \textsl{(Rademacher) cotype $q$} 
if there is a constant $C_q$ such that for all finite sequences $\{x_i\}$ in $\B$,
$$
	\Ex \Bigl\| \sum_{i=1}^n \epsilon_i x_i \Bigr\|^q
	\ \ge \ C_q^{-q} \left( \sum_{i=1}^n  \| x_i \|^q \right)^{1/q} .
$$
Ledoux and Talagrand \cite{LedouxTalagrand:91}, page 247, note that type 
and cotype properties appear as dual notions: If a Banach space $\B$ is of type $p$, 
its dual space $\B'$ is of cotype $q = p/(p-1)$.   

One of the basic results concerning Banach spaces with type $p$ and 
cotype $q$ is the following proposition: 

\begin{proposition}
\label{prop: Hoffmann-Jorgensen}
\cite[Proposition 9.11, page 248]{LedouxTalagrand:91}.

\noindent
If $\B$ is of type $p \ge 1$ with constant $T_p$, then
$$
	\Ex \|S_n\|^p \ \le \ (2T_p)^p \sum_{i=1}^n \Ex \|X_i\|^p .
$$
If $\B$ is of cotype $q \ge 1$ with constant $C_q$, then
$$
	\Ex \|S_n\|^q \ \ge \ (2C_q)^{-q} \sum_{i=1}^n \Ex \|X_i\|^q .
$$
\end{proposition}

As shown in \cite{LedouxTalagrand:91}, page 27, the Banach space $L_r(\mu)$ 
with $1 \le r < \infty$ (cf.\ section~\ref{subsec: L_r-spaces}) is of type $\min(r,2)$. 
Similarly, $L_r(\mu)$ is co-type $\max(r,2)$. In case of $r \ge 2 = p$, explicit values 
for the constant $T_p$ in Proposition~\ref{prop: Hoffmann-Jorgensen} can be 
obtained from the optimal constants in Khintchine's inequalities due to \cite{Haagerup:81}.

\begin{lemma}
\label{lem: Haagerup}
For $2 \le r < \infty$, the space $L_r(\mu)$ is of type $2$ with constant $T_2 = B_r$, where
$$
	B_r \ := \ 2^{1/2} \left ( \frac{\Gamma ((r+1)/2)}{\sqrt{\pi}} \right )^{1/r} .
$$
\end{lemma}

\begin{corollary}
\label{cor: Haagerup}
For $\B = L_r(\mu)$, $2 \le r < \infty$, inequality~(\ref{ineq: Nemirovski}) is satisfied with $K = 4 B_r^2$.
\end{corollary}

Note that $B_2 = 1$ and 
$$
	\frac{B_r}{\sqrt{r}} \ \to \ \frac{1}{\sqrt{e}}
	\quad \mbox{as} \ r \to \infty .
$$
Thus for large values $r$, the conclusion of Corollary~\ref{cor: Haagerup} is 
weaker than the one of Corollary~\ref{cor: Nemirovski 1}.

\subsection{The space $\ell_\infty^d$}

The preceding results apply only to $r < \infty$, so the special space 
$\ell_\infty^d$ requires different arguments. At first we deduce a new type 
inequality based on Hoeffding's \cite{MR0144363} exponential inequality: 
If $\epsilon_1, \epsilon_2, \ldots, \epsilon_n$ are independent Rademacher 
random variables, $a_1, a_2, \ldots , a_n$ are real numbers and 
$v^2 := \sum_{i=1}^n a_i^2$, then the tail probabilities of the random 
variable $\Bigl| \sum_{i=1}^n a_i \epsilon_i \Bigr|$ may be bounded as follows:
\begin{equation}
\label{eq:Hoeffding}
	\Pr \left( \Bigl| \sum_{i=1}^n a_i \epsilon_i \Bigr| \ge z \right)
	\ \le \ 2 \exp \Bigl( - \frac{z^2}{2 v^2} \Bigr) ,  \quad z \ge 0 .
\end{equation}
At the heart of these tail bounds is the following exponential moment bound:
\begin{equation}
\label{eq:Hoeffding2}
	\Ex \exp \Bigl( t \sum_{i=1}^n a_i \epsilon_i \Bigr)
	\ \le \ \exp(t^2 v^2 / 2) ,  \quad t \in \R .
\end{equation}
From the latter bound we shall deduce the following type inequality in Section~\ref{sec: Proofs}:

\begin{lemma}
\label{lem: Veraar}
The space $\ell_\infty^d$ is of type $2$ with constant $\sqrt{2 \log(2d)}$.
\end{lemma}

Using this upper bound together with Proposition~\ref{prop: Hoffmann-Jorgensen} 
yields another Nemirovski type inequality:

\begin{corollary}
\label{cor: Veraar}
For $(\B,\|\cdot\|) = \ell_\infty^d$, inequality~(\ref{ineq: Nemirovski}) 
is satisfied with $K = \KT2(d,\infty) = 8 \log(2d)$.
\end{corollary}

\paragraph{Refinements.}
Let $T_2(\ell_\infty^d)$ be the optimal type 2 constant for the 
space $\ell_\infty^d$. So far we know that $T_2(\ell_\infty^d) \le \sqrt{2 \log(2d)}$. 
Moreover, by a modification of Example~\ref{ex: lower bound ell_infty} one can show that
\begin{equation}
\label{ineq: type2 ge cd}
	T_2(\ell_\infty^d) \ \ge \ c_d \ := \ \sqrt{ \Ex \max_{1 \le j \le d} Z_j^2 } .
\end{equation}

The constants $c_d$ can be expressed or bounded in terms of the 
distribution function $\Phi$ of $N(0,1)$, i.e. $\Phi(z) = \int_{-\infty}^z \phi(x) \, dx$ 
with $\phi(x) = \exp(- x^2/2) / \sqrt{2\pi}$. Namely, with $W := \max_{1 \le j \le d} |Z_j|$,
$$
	c_d^2 \ = \ \Ex(W^2)
	\ = \ \Ex \int_0^\infty 2t 1_{[t \le W]} \, dt
	\ = \ \int_0^\infty 2t \Pr(W \ge t) \, dt ,
$$
and for any $t > 0$,
$$
	\Pr(W \ge t) \ \begin{cases}
		= \ 1 - \Pr(W < t) \ = \ 1 - \Pr(|Z_1| < t)^d \ = \ 1 - (2 \Phi(t) - 1)^d , \\
		\le \ d \Pr(|Z_1| \ge t) \ = \ 2d (1 - \Phi(t)) .
	\end{cases}
$$
These considerations and various bounds for $\Phi$ will allow us to derive explicit bounds for $c_d$. 

On the other hand, Hoeffding's inequality (\ref{eq:Hoeffding}) 
has been refined by Pinelis \cite{MR1272088,MR2339301} as follows:
\begin{equation}
\label{eq:Pinelis}
	\Pr \left( \Bigl| \sum_{i=1}^n a_i \epsilon_i \Bigr| \ge z \right)
	\ \le \ 2K (1 - \Phi(z/v)), \quad z > 0 , 
\end{equation}
where $K$ satisfies $3.18 \le K \le 3.22$. This will be the main ingredient 
for refined upper bounds for $T_2(\ell_\infty^d)$. 
The next lemma summarizes our findings:

\begin{lemma}
\label{lem: Jon Mark}
The constants $c_d$ and $T_2(\ell_\infty^d)$ satisfy the following inequalities:
\begin{eqnarray}
\label{RademacherTypeTwoBestBounds}
	\sqrt{2 \log d + h_1(d)} \le c_d \le 
	\left \{ \begin{array}{l l} T_2(\ell_{\infty}^d ) \le \sqrt{2 \log d + h_2(d)}, & \qquad d \ge 1 \\
	                                     \sqrt{2\log d}, & \qquad d \ge 3 \\
	                                     \sqrt{2\log d + h_3 (d)}, &\qquad d\ge 1 
	         \end{array} \right .	                                   
\end{eqnarray}
where $h_2(d) \le 3$, $h_2(d)$ becomes negative for 
$d > 4.13795 \times 10^{10}$, 
$h_3(d)$ becomes negative for $d\ge 14$, 
and $h_j (d) \sim - \log \log d$ as $d\rightarrow \infty$ for $j=1,2,3$.
\end{lemma}

In particular, one could replace $\KT2(d,\infty)$ in Corollary~\ref{cor: Veraar} with $8 \log d + 4 h_2(d)$.

\section{The empirical process approach:
	Truncation and\\ Bernstein's inequality.}
\label{sec: Bernstein}

An alternative to Hoeffding's exponential tail inequality (\ref{eq:Hoeffding}) 
is a classical exponential bound due to Bernstein (see e.g. \cite{Benn:prob:1962}): 
Let $Y_1, Y_2, \ldots, Y_n$ be independent random variables with mean zero 
such that $|Y_i| \le \kappa$. Then for any $v^2 \ge \sum_{i=1}^n \Var(Y_i)$,
\begin{eqnarray}
	\Pr \left( \Bigl| \sum_{i=1}^n Y_i \Bigr| \ge  x \right)
	\le 2 \exp \left ( - \frac{x^2}{2 ( v^2 + \kappa x/3)} \right ) , \qquad x > 0 .
\label{BernsteinInequality}
\end{eqnarray}
We will not use this inequality itself but rather an exponential 
moment inequality underlying its proof:

\begin{lemma}
\label{lem: linexp}
For $L > $ define
$$
	{\rm e}(L) \ := \ \exp(1/L) - 1 - 1/L .
$$
Let $Y$ be a random variable with mean zero and variance 
$\sigma^2$ such that $|Y| \le \kappa$. Then for any $L>0$,
$$
	\Ex \exp \Bigl( \frac{Z}{\kappa L} \Bigr)
	\ \le \ 1 + \frac{\sigma^2 {\rm e}(L)}{\kappa^2}
	\ \le \ \exp \Bigl( \frac{\sigma^2 {\rm e}(L)}{\kappa^2} \Bigr) .
$$
\end{lemma}

With the latter exponential moment bound we can prove a 
moment inequality for random vectors with bounded components:

\begin{lemma}
\label{lem: Bernstein}
Suppose that $X_i = (X_{i,j})_{j=1}^d$ satisfies $\|X_i\|_\infty \le \kappa$, 
and let $\Gamma$ be an upper bound for 
$\max_{1 \le j \le d} \sum_{i=1}^n \Var(X_{i,j})$. Then for any $L>0$,
$$
	\sqrt{ \Ex \|S_n\|_\infty^2 }
	\ \le \ \kappa L \log (2d) + \frac{\Gamma L  \, {\rm e}(L)}{\kappa} .
$$
\end{lemma}

Now we return to our general random vectors $X_i \in \R^d$ with 
mean zero and $\Ex \|X_i\|_\infty^2 < \infty$. 
They are split into two random vectors via truncation: $X_i = X_i^{(a)} + X_i^{(b)}$ with
$$
	X_i^{(a)} \ := \ 1_{[\|X_i\|_\infty \le \kappa_o]} X_i
	\quad\text{and}\quad
	X_i^{(b)} \ := \ 1_{[\|X_i\|_\infty > \kappa_o]} X_i
$$
for some constant $\kappa_o > 0$ to be specified later. 
Then we write $S_n = A_n + B_n$ with the centered random sums
$$
	A_n \ := \ \sum_{i=1}^n (X_i^{(a)} - \Ex X_i^{(a)})
	\quad\text{and}\quad
	B_n \ := \ \sum_{i=1}^n (X_i^{(b)} - \Ex X_i^{(b)}) .
$$
The sum $A_n$ involves centered random vectors in $[-2\kappa_o, 2\kappa_o]^d$ 
and will be treated by means of Lemma~\ref{lem: Bernstein}, while $B_n$ will be 
bounded with elementary methods. Choosing the threshold $\kappa$ 
and the parameter $L$ carefully yields the following theorem.

\begin{theorem}
\label{thm: Bernstein-Nemirovski}
In case of $(\B,\|\cdot\|) = \ell_\infty^d$ for some $d \ge 1$, 
inequality~(\ref{ineq: Nemirovski}) holds with
$$
	K \ = \KTB(d,\infty) := \bigl( 1 + 3.46 \sqrt{\log(2d)} \bigr)^2 .
$$
If the random vectors $X_i$ are symmetrically distributed around $0$, one may even set
$$
	K \ = \ \KTB^{(\rm symm)}(d,\infty) \ = \ \bigl( 1 + 2.9 \sqrt{\log(2d)} \bigr)^2 .
$$
\end{theorem}

\section{Comparisons.}
\label{sec: Comparisons}

In this section we compare the three approaches just described for the 
space $\ell_\infty^d$. As to the random vectors $X_i$, we broaden our 
point of view and consider three different cases:
\begin{description}
\item[General case:] The random vectors $X_i$ are independent with 
$\Ex \|X_i\|_\infty^2 < \infty$ for all $i$.
\item[Centered case:] In addition, $\Ex X_i = 0$ for all $i$.
\item[Symmetric case:] In addition, $X_i$ is symmetrically distributed around $0$ for all $i$.
\end{description}
In view of the general case, we reformulate inequality~(\ref{ineq: Nemirovski}) as follows:
\begin{equation}
\label{ineq: Nemirovski general}
	\Ex \| S_n - \Ex S_n \|_\infty^2
	\ \le \ K \sum_{i=1}^n \Ex \|X_i\|_\infty^2 .
\end{equation}
One reason for this extension is that in some applications, particularly in 
connection with empirical processes, it is easier and more natural to work 
with uncentered summands $X_i$. Let us discuss briefly the consequences 
of this extension in the three frameworks:

\paragraph{Nemirovski's approach:}
Between the centered and symmetric case there is no difference. If (\ref{ineq: Nemirovski}) 
holds in the centered case for some $K$, then in the general case
$$
	\Ex \| S_n - \Ex S_n \|_\infty^2
	\ \le \ K \sum_{i=1}^n \Ex \|X_i - \Ex X_i\|_\infty^2
	\ \le \ 4K \sum_{i=1}^n \Ex \|X_i\|_\infty^2 .
$$
The latter inequality follows from the general fact that
$$
	\Ex \|Y - \Ex Y\|^2 \ \le \ \Ex \bigl( (\|Y\| + \|\Ex Y\|)^2 \bigr)
	\ \le \ 2 \Ex \|Y\|^2 + 2 \|\Ex Y\|^2 \ \le \ 4 \Ex \|Y\|^2 .
$$
This looks rather crude at first glance, but in case of the maximum norm 
and high dimension $d$, the factor $4$ cannot be reduced. For let $Y \in \R^d$ 
have independent components $Y_1, \ldots, Y_d \in \{-1,1\}$ with
$\Pr(Y_j = 1) = 1 - \Pr(Y_j = -1) = p \in [1/2, 1)$. Then $\|Y\|_\infty \equiv 1$, 
while $\Ex Y = (2p - 1)_{j=1}^d$ and
$$
	\|Y - \Ex Y\|_\infty \ = \ \begin{cases}
		2(1 - p) & \text{if} \ Y_1 = \cdots = Y_d = 1 , \\
		2p       & \text{else} .
	\end{cases}
$$
Hence
$$
	\frac{\Ex \|Y - \Ex Y\|_\infty^2}{\Ex \|Y\|_\infty^2}
	\ = \ 4 \bigl( (1 - p)^2 p^d + p^2 (1 - p^d) \bigr) .
$$
If we set $p = 1 - d^{-1/2}$ for $d \ge 4$, then the latter ratio converges to $4$ as $d \to \infty$.

\paragraph{The approach via Rademacher type 2 inequalities:}
The first part of Proposition~\ref{prop: Hoffmann-Jorgensen}, involving the 
Rademacher type constant $T_p$, remains valid if we drop the assumption 
that $\Ex X_i = 0$ and replace $S_n$ with $S_n - \Ex S_n$. Thus there is 
no difference between the general and the centered case. In the symmetric case, 
however, the factor $2^p$ in Proposition~\ref{prop: Hoffmann-Jorgensen} becomes 
superfluous. Thus, if (\ref{ineq: Nemirovski}) holds with a certain constant $K$ in the 
general and centered case, we may replace $K$ with $K/4$ in the symmetric case.

\paragraph{The approach via truncation and Bernstein's inequality:}
Our proof for the centered case does not utilize that $\Ex X_i = 0$, so again there is 
no difference between the centered and general case. However, in the symmetric case, 
the truncated random vectors $1\{\|X_i\|_\infty \le \kappa\} X_i$ and 
$1\{\|X_i\|_\infty > \kappa\} X_i$ are centered, too, which leads to the substantially smaller 
constant $K$ in Theorem~\ref{thm: Bernstein-Nemirovski}.

\begin{figure}[h]
\centering
\includegraphics[width=0.7\textwidth]{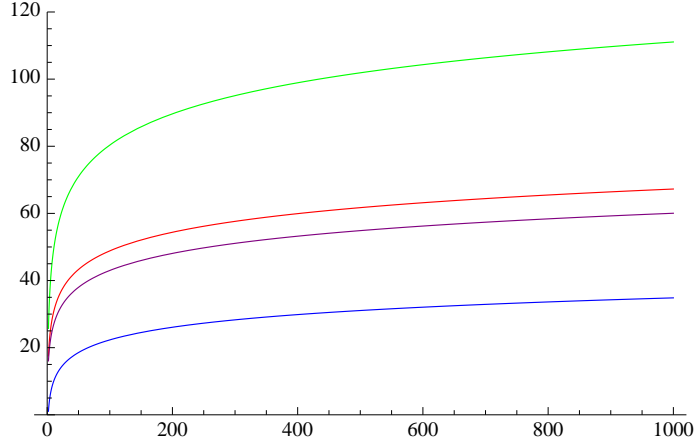}
\caption{Comparison of $K(d, \infty)$ obtained via the three proof methods:  
   Blue (bottom) $=$ Nemirovski; Magenta and Red (middle) $=$ type 2 inequalities; 
   Green (top) $ = $ truncation and Bernstein inequality}
\label{fig: Comparisons}
\end{figure}

\paragraph{Summaries and comparisons.}
Table~\ref{tab: All Ks 1} summarizes the constants $K = K(d,\infty)$ we have found 
so far by the three different methods and for the three different cases. 
Table~\ref{tab: All Ks 2} contains the corresponding limits
$$
	K^* \ := \ \lim_{d \to \infty} \, \frac{K(d,\infty)}{\log d} .
$$
Interestingly, there is no global winner among the three methods. 
But for the centered case, Nemirovski's approach yields asymptotically 
the smallest constants. In particular,
\begin{eqnarray*}
	\lim_{d\to\infty} \, \frac{\KTB(d,\infty)}{\KNem(d,\infty)} 
	& = & \frac{3.46^2}{2e} \deq 2.20205 , \\
	\lim_{d\to\infty} \, \frac{\KT2(d,\infty)}{\KNem(d,\infty)}
	& = & \frac{4}{e} \deq 1.47152 , \\
	\lim_{d\to\infty} \, \frac{\KTB(d,\infty)}{\KT2(d,\infty)}
	& = & \frac{3.46^2}{8} \deq 1.49645 .
\end{eqnarray*}
The conclusion at this point seems to be that Nemirovski's approach and the 
type 2 inequalities yield better constants than Bernstein's inequality and truncation. 
Figure~\ref{fig: Comparisons} shows the constants $K(d,\infty)$ for the centered 
case over a certain range of dimensions $d$.

\begin{table}[h]
\centering
\begin{tabular}{|c||c|c|c|}
	\hline
		& General case
			& Centered case
				& Symmetric case \\
	\hline\hline
	Nemirovski
		& $8e \log d - 4e$
			& $2e \log d - e$
				& $2e \log d - e$  \\
	\hline
	Type 2 inequalities
		& $8 \log(2d)$
			& $8 \log(2d)$
				& $ 2 \log(2d)$ \\
		& $8 \log d + 4 h_2(d)$
			& $8 \log d + 4 h_2(d)$
				& $ 2 \log d + h_2(d)$ \\
	\hline
	Truncation/Bernstein
		& \!\!\!$\bigl( 1 + 3.46 \sqrt{\log(2d)} \bigr)^{\!2}$\!\!
			& \!\!\!$\bigl( 1 + 3.46 \sqrt{\log(2d)} \bigr)^{\!2}$\!\!
		    		& \!\!\!$\bigl( 1 + 2.9 \sqrt{\log(2d)} \bigr)^{\!2}$\!\! \\
	\hline
\end{tabular}
\caption{The different constants $K(d,\infty)$.}
\label{tab: All Ks 1}
\end{table}

\begin{table}[h]
\centering
\begin{tabular}{|c||c|c|c|}
	\hline
		& General case
			& Centered case
				& Symmetric case \\
	\hline\hline
	Nemirovski
		& $8e \deq 21.7463$
			& $2e \deq \textbf{5.4366}$
				& $2e \deq 5.4366$ \\
	\hline
	Type 2 inequalities
		& $\textbf{8.0}$
			& $8.0$
				& $\textbf{2.0}$ \\
	\hline
	Truncation/Bernstein
		& $3.46^2 = 11.9716$
			& $3.46^2 = 11.9716$
				& $2.9^2 = 8.41$ \\
	\hline
\end{tabular}
\caption{The different limits $K^*$.}
\label{tab: All Ks 2}
\end{table}

\section{Proofs.}
\label{sec: Proofs}

\subsection{Proofs for Section~\ref{sec: Geometry}}

\paragraph{Proof of (\ref{ineq: norms}).}

In case of $r = \infty$, the asserted inequalities read
$$
	\|x\|_\infty \ \le \ \|x\|_q \ \le \ d^{1/q} \|x\|_\infty
	\quad\text{for} \ 1 \le q < \infty
$$
and are rather obvious. For $1 \le q < r < \infty$, (\ref{ineq: norms}) is an easy 
consequence of H\"older's inequality. \hfill $\Box$

\paragraph{Proof of Lemma~\ref{lem: good news!}.}
In case of $r = 2$, $V(f + g)$ is equal to $V(f) + DV(f,g) + V(g)$. 
In case of $r \ge 2$ and $\|f\|_r = 0$, both $DV(f,g)$ and $\int h(f) g \, d\mu$ are 
equal to zero, and the asserted inequalities reduce to the trivial statement 
that $V(g) \le (r - 1) V(g)$. Thus let us restrict our attention to the case $r > 2$ and $\|f\|_r > 0$.

Note first that the mapping
$$
	\R \ni t \ \mapsto \ h_t := |f + tg|^r
$$
is pointwise twice continuously differentiable with derivatives
\begin{eqnarray*}
	\dot{h}_t  & = & r |f + tg|^{r-1} \mathrm{sign}(f + tg) g
		\ = \ r |f + tg|^{r-2} (f + tg) g , \\
	\ddot{h}_t & = & r(r-1) |f + tg|^{r-2} g^2 .
\end{eqnarray*}
By means of the inequality $|x + y|^b \le 2^{b-1} \bigl( |x|^b + |y|^b \bigr)$ 
for real numbers $x$, $y$ and $b \ge 1$, a consequence of Jensen's 
inequality, we can conclude that for any bound $t_o > 0$,
\begin{eqnarray*}
	\max_{|t| \le t_o} |\dot{h}_t|
	& \le & r 2^{r-2} \bigl( |f|^{r-1} |g| + t_o^{r-1}|g|^r \bigr) , \\
	\max_{|t| \le t_o} |\ddot{h}_t|
	& \le & r(r-1) 2^{r-3} \bigl( |f|^{r-2} |g|^2 + t_o^{r-2}|g|^r \bigr) .
\end{eqnarray*}
The latter two envelope functions belong to $L_1(\mu)$. 
This follows from H\"older's inequality which we rephrase for our purposes in the form
\begin{equation}
\label{eq:Hoelder}
	\int |f|^{(1 - \lambda) r} |g|^{\lambda r} \, d\mu
	\ \le \ \|f\|_r^{(1 - \lambda) r} \|g\|_r^{\lambda r}
	\quad\text{for} \ 0 \le \lambda \le 1 .
\end{equation}
Hence we may conclude via dominated convergence that
$$
	t \ \mapsto \ \tilde{v}(t) := \|f + tg\|_r^r
$$
is twice continuously differentiable with derivatives
\begin{eqnarray*}
	\tilde{v}'(t)  & = & r \int |f + tg|^{r-2} (f + tg) g \, d\mu , \\
	\tilde{v}''(t) & = & r (r-1) \int |f + tg|^{r-2} g^2 \, d\mu .
\end{eqnarray*}
This entails that
$$
	t \ \mapsto \ v(t) := V(f + tg) \ = \ \tilde{v}(t)^{2/r}
$$
is continuously differentiable with derivative
$$
	v'(t) \ = \ (2/r) \tilde{v}(t)^{2/r - 1} \tilde{v}'(t)
	\ = \ \tilde{v}^{2/r-1} (t) \int h(f + tg) g \, d\mu .
$$
For $t = 0$ this entails the asserted expression for $DV(f,g)$. 
Moreover, $v(t)$ is twice continuously differentiable on the set 
$\{t \in \R : \|f + tg\|_r > 0\}$ which equals either $\R$ or 
$\R \setminus \{t_o\}$ for some $t_o \ne 0$. On this set the second derivative equals
\begin{eqnarray*}
	v''(t)
	& = & (2/r) \tilde{v}(t)^{2/r - 1} \tilde{v}''(t)
		+ (2/r)(2/r - 1) \tilde{v}(t)^{2/r - 2} \tilde{v}'(t)^2 \\
	& = & 2 (r - 1) \int \frac{|f + tg|^{r-2}}{\|f + tg\|_r^{r-2}} g^2 \, d\mu
		- 2 (r - 2)
			\Bigl( \int \frac{|f + tg|^{r-2} (f + tg)}{\|f + tg\|_r^{r-1}} g \, d\mu \Bigr)^2 \\
	& \le & 2 (r - 1) \int \Bigl| \frac{f + tg}{\|f + tg\|_r} \Bigr|^{r - 2} |g|^2 \, d\mu \\
	& \le & 2 (r - 1) \|g\|_r^2 \ = \ 2 (r - 1) V(g)
\end{eqnarray*}
by virtue of H\"older's inequality \eqref{eq:Hoelder} with $\lambda = 2/r$. Consequently,
by using 
\begin{eqnarray*}
    v'(t) - v'(0) = \int_0^t v'' (s) \, ds \le 2 (r-1) V(g) t ,
\end{eqnarray*}
we find that 
\begin{eqnarray*}
\lefteqn{V(f + g) - V(f) - DV(f,g)} \\
	& = & v(1) - v(0) - v'(0) 
	          =  \int_0^1 (v'(t) - v'(0)) \, dt \\
	& \le & 2 (r - 1) V(g) \int_0^1 t \, dt 
	           =  (r - 1) V(g) .
\end{eqnarray*}\\[-6ex]
\strut\hfill	$\Box$

\paragraph{Proof of Theorem~\ref{thm: Nemirovski}.}
The first part is an immediate consequence of the considerations 
preceding the theorem. It remains to prove the (in)equalities and expansion 
for $\KNem(d,r)$. Note that $\KNem(d,r)$ is the infimum of $h(q) d^{-2/r}$ 
over all real $q \in [2,r]$, where $h(q) := (q - 1) d^{2/q}$ satisfies the equation
$$
	h'(q)
	\ = \ \frac{d^{2/q}}{q^2} \bigl( (q - \log d)^2 - (\log d - 2) \log d \bigr) .
$$
Since $7 < e^2 < 8$, this shows that $h$ is strictly increasing on $[2,\infty)$ if $d \le 7$. Hence
$$
	\KNem(d,r) \ = \ h(2) d^{-2/r} \ = \ d^{1 - 2/r}	\quad\text{if} \ d \le 7 .
$$
For $d \ge 8$, one can easily show that $\log d - \sqrt{(\log d - 2)\log d} < 2$, 
so that $h$ is strictly decreasing on $[2, r_d]$ and strictly increasing on $[r_d, \infty)$, where
$$
	r_d \ := \ \log d + \sqrt{(\log d - 2)\log d}
	\ \begin{cases}
		< \ 2 \log d , \\
		> \ 2 \log d - 2 .
	\end{cases}
$$
Thus for $d \ge 8$,
$$
	\KNem(d,r) \ = \ \begin{cases}
		h(r) d^{-2/r} \ = \ r - 1 \ < \ 2 \log d - 1
			& \text{if} \ r \le r_d , \\
		h(r_d) d^{-2/r} \ \le \ h(2 \log d) \ = \ 2e \log d - e
			& \text{if} \ r \ge r_d .
	\end{cases}
$$
Moreover, one can verify numerically that $\KNem(d,r) \le d \le 2e \log d - e$ for $3 \le d \le 7$.

Finally, for $d \ge 8$, the inequalities $r_d' := 2 \log d - 2 < r_d < r_d'' := 2 \log d$ yield
$$
	\KNem(d,\infty) = h(r_d)
	\ \ge \ (r_d' - 1) d^{2/r_d''}
	\ = \ 2e \log d - 3e ,
$$
and for $1 \le d \le 7$, the inequality $d = \KNem(d,\infty) \ge 2e \log(d) - 3e$ 
is easily verified.\hfill$\Box$

\subsection{Proofs for Section~\ref{sec: Type and cotype}}

\paragraph{Proof of Lemma~\ref{lem: Haagerup}.}
The following proof is standard; see e.g. 
\cite{MR576407}, page 160, \cite{LedouxTalagrand:91}, page 247.
Let $x_1, \ldots , x_n$ be fixed functions in $L_r(\mu)$. 
Then by \cite{Haagerup:81}, for any $t \in T$,
\begin{equation}
\label{KhinchineInequalityWithOptimalConst}
	\biggl\{ \Ex \Bigl| \sum_{i=1}^n \epsilon_i x_i(t) \Bigr|^r \biggr\}^{1/r}
	\ \le \ B_r \Bigl( \sum_{i=1}^n |x_i(t)|^2 \Bigr)^{1/2} .
\end{equation}
To use inequality (\ref{KhinchineInequalityWithOptimalConst}) for finding 
an upper bound for the type constant for $L_r$, rewrite it as
$$
	\Ex \Bigl| \sum_{i=1}^n \epsilon_i x_i(t) \Bigr|^r
	\ \le \ B_r^r \Bigl( \sum_{i=1}^n |x_i(t)|^2 \Bigr)^{r/2} .
$$
It follows from Fubini's theorem and the previous inequality that
\begin{eqnarray*}
	\Ex \Bigl\| \sum_{i=1}^n \epsilon_i x_i \Bigr\|_r^r 
	& = & \Ex \int \Bigl| \sum_{i=1}^n \epsilon_i x_i(t) \Bigr|^r \, d\mu(t) \\
	& = & \int \Ex \Bigl| \sum_{i=1}^n \epsilon_i x_i (t)|^r \, d\mu(t) \\
	& \le & B_r^r \int \biggl( \sum_{i=1}^n |x_i(t)|^2 \biggr)^{r/2} \, d\mu(t) .
\end{eqnarray*}
Using the triangle inequality (or Minkowski's inequality), we obtain
\begin{eqnarray*}
	\biggl\{ \Ex \Bigl\| \sum_{i=1}^n \epsilon_i x_i \Bigr\|_r^r \biggr\}^{2/r}  
	& \le & B_r^2 \biggl\{ \int \Bigl( \sum_{i=1}^n |x_i(t)|^2 \Bigr)^{r/2} 
		\, d\mu(t) \biggr\}^{2/r} \\
	& \le & B_r^2 \sum_{i=1}^n \Bigl( \int |x_i(t)|^r \, d\mu(t) \Bigr)^{2/r} \\
	& = & B_r^2 \sum_{i=1}^n \|x_i\|_r^2 .
\end{eqnarray*}
Furthermore, since $g(v) = v^{2/r}$ is a concave function of $v \ge 0$, 
the last display implies that
$$
	\Ex \Bigl\| \sum_{i=1}^n \epsilon_i x_i \Bigr\|_r^2
	\ \le \ \biggl\{ \Ex \Bigl\| \sum_{i=1}^n \epsilon_i x_i \Bigr\|_r^r \biggr\}^{2/r}   
	\ \le \ B_r^2 \sum_{i=1}^n \|x_i\|_r^2 .
	\eqno{\Box}
$$

\paragraph{Proof of Lemma~\ref{lem: Veraar}.}
For $1 \le i \le n$ let $x_i = (x_{im})_{m=1}^d$ be an arbitrary fixed vector in 
$\R^d$, and set $S := \sum_{i=1}^n \epsilon_i x_i$. Further let $S_m$ be the 
$m$-th component of $S$ with variance $v_m^2 := \sum_{i=1}^n x_{im}^2$, 
and define $v^2 := \max_{1 \le m \le d} v_m^2$, which is not greater than 
$\sum_{i=1}^n \|x_i\|_\infty^2$. It suffices to show that
$$
	\Ex \|S\|_\infty^2 \ \le \ 2 \log(2d) v^2 .
$$
To this end note first that $h : [0,\infty) \to [1,\infty)$ with 
$h(t) := \cosh(t^{1/2}) = \sum_{k=0}^\infty t^k / (2k)!$ is bijective, increasing and convex. 
Hence its inverse function $h^{-1} : [1,\infty) \to [0,\infty)$ is increasing and concave, 
and one easily verifies that $h^{-1}(s) = \bigl( \log(s + (s^2 - 1)^{1/2}) \bigr)^2 \le (\log(2s))^2$. 
Thus it follows from Jensen's inequality that for arbitrary $t > 0$,
\begin{eqnarray*}
	\Ex \|S\|_\infty^2
	& = & t^{-2} \Ex h^{-1} \bigl( \cosh(\|tS\|_\infty) \bigr)
		\ \le \ t^{-2} h^{-1} \bigl( \Ex \cosh(\|tS\|_\infty) \bigr) \\
	& \le & t^{-2} \bigl( \log \bigl( 2 \Ex \cosh(\|tS\|_\infty) \bigr) \bigr)^2 .
\end{eqnarray*}
Moreover,
$$
	\Ex \cosh(\|tS\|_\infty)
	\ = \ \Ex \max_{1 \le m \le d} \cosh(t S_m)
	\ \le \ \sum_{m=1}^d \Ex \cosh(t S_m)
	\ \le \ d \exp(t^2 v^2/2) ,	
$$
according to (\ref{eq:Hoeffding2}), whence
$$
	\Ex \|S\|_\infty^2
	\ \le \ t^{-2} \log \bigl( 2d \exp(t^2 v^2/2) \bigr)^2
	\ = \ \bigl( \log(2d)/t + t v^2 / 2 \bigr)^2 .
$$
Now the assertion follows if we set $t = \sqrt{2\log(2d)/v^2}$.	\hfill	$\Box$

\paragraph{Proof of (\ref{ineq: type2 ge cd}).}
We may replace the random sequence $\{X_i\}$ in 
Example~\ref{ex: lower bound ell_infty} with the random sequence 
$\{\epsilon_i X_i\}$, where $\{\epsilon_i\}$ is a Rademacher 
sequence independent of $\{X_i\}$. Thereafter we condition on $\{X_i\}$, i.e.\ 
we view it as a deterministic sequence such that $n^{-1} \sum_{i=1}^n X_i^{}X_i^\top$ 
converges to the identity matrix $I_d$ as $n \to \infty$, by the strong law of large numbers. 
Now Lindeberg's version of the multivariate Central Limit Theorem shows that
$$
	\sup_{n \ge 1} \frac{\Ex \Bigl\| \sum_{i=1}^n \epsilon_i X_i \Bigr\|_\infty^2}
	                              {\sum_{i=1}^n \|X_i\|_\infty^2}
	\ \ge \ \sup_{n \ge 1} \Ex \Bigl\| n^{-1/2} \sum_{i=1}^n \epsilon_i X_i \Bigr\|_\infty^2
	\ \ge \ c_d^2 .
	\eqno{\Box}
$$

\paragraph{Inequalities for $\Phi$.}
The subsequent results will rely on (\ref{eq:Pinelis}) and several inequalities 
for $1 - \Phi(z)$. The first of these is:
\begin{eqnarray}
\label{MillsRatio}
	1 - \Phi(z) \le z^{-1} \phi (z) ,	\qquad z > 0, 
\end{eqnarray}
which is known as {\sl Mills' ratio}; see \cite{MR0005558} and \cite{MR1906389} 
for related results. The proof of this upper bound is easy: 
Since $\phi'(z) = - z \phi (z)$ it follows that
\begin{equation}
\label{MillsRatioUpperBoundProof}
	1 - \Phi(z)
	= \int_z^{\infty} \phi (t) \, dt
	\le \int_z^{\infty} \frac{t}{z} \phi(t) \, dt
	= \frac{-1}{z}\int_z^\infty \phi'(t) \, dt
	= \frac{\phi(z)}{z} .
\end{equation}
A very useful pair of upper and lower bounds for $1-\Phi(z)$ are as follows:
\begin{eqnarray}
\label{ineq:KMSW}
	\frac{2}{z + \sqrt{z^2+4}} \phi (z)
	\le 1 - \Phi(z)
	\le \frac{4}{3z + \sqrt{z^2 +8}} \phi (z) , \qquad z > -1 ;
\end{eqnarray}
the inequality on the left is due to Komatsu (see e.g. 
\cite{MR0345224} p.\ 17), while the inequality on the right is 
an improvement of an earlier result of Komatsu due to \cite{MR1677442}.

\paragraph{Proof of Lemma~\ref{lem: Jon Mark}.}
%
To prove the upper bound for $T_2(\ell_\infty^d)$, let $(\epsilon_i)_{i\geq 1}$ be a Rademacher sequence. 
With $S$ and $S_m$ as in the proof of Lemma~\ref{lem: Veraar}, we may write
\begin{eqnarray*}
	\Ex \|S\|_\infty^2
	&=& \int_0^\infty 2t
                     \Pr\Big(\sup_{1\leq m\leq d} |S_m|>t \Big) \, dt \\ 
	&\le& \delta^2 + \int_\delta^\infty 2t
            \Pr\Big(\sup_{1\leq m\leq d} |S_m| >t \Big) \, dt \\ 
	&\le& \delta^2 + \sum_{m=1}^d \int_\delta^\infty 2t
            \Pr\big( |S_m| >t \big) \, dt.
\end{eqnarray*}
Now by \eqref{eq:Pinelis} 
with $v^2$ and $v_m^2$ as in the proof of Lemma~\ref{lem: Veraar}, 
followed by Mills' ratio \eqref{MillsRatio},
\begin{eqnarray}
	\int_\delta^{\infty} 2 t \Pr( | S_m | > t ) \, dt
	& \le & \int_{\delta}^{\infty} \frac{4K v_m}{ \sqrt{2 \pi} t} t e^{- t^2 /(2v_m^2)} \, dt 
		\nonumber \\
	& = &  \frac{4K v_m }{ \sqrt{2 \pi} } \int_{\delta}^{\infty} e^{- t^2/(2v_m^2)} \, dt 
		\ = \ 4K v_m^2 \int_{\delta}^{\infty} \frac{e^{-t^2/(2v_m^2)}}{\sqrt{2\pi} v_m} \, dt 
		\nonumber \\
	& = & 4K v_m^2  (1 - \Phi (\delta / v_m))
		\ \le \ 4 K v^2 (1 - \Phi (\delta /v)) .
	\label{ThirdTryTailOfNormalBnd}
\end{eqnarray}
Now instead of the Mills' ratio bound (\ref{MillsRatio}) for the tail of the normal distribution, we use
the upper bound part of (\ref{ineq:KMSW}) due to \cite{MR1677442}. This yields
\[
	\int_\delta^\infty 2t \Pr(|S_m| >t)\, dt
	\le 4K v^2 (1-\Phi(\delta/v))
	\le \frac{4c v^2}{3 \delta/v + \sqrt{\delta^2/v^2 + 8} } e^{-\delta^2/(2v^2)},
\]
where we have defined $c := 4 K / \sqrt{2\pi} = 12.88/ \sqrt{2 \pi}$, and hence
\[
	\Ex \|S\|^2\leq \delta^2 + \frac{4 c d v^2}{3 \delta/v + \sqrt{\delta^2/v^2 + 8}}
		e^{-\delta^2/(2v^2)} .
\]
Taking 
\[
	\delta^2 = v^2 2 \log \left ( \frac{cd/2}{\sqrt{2 \log (cd/2)}} \right )
\]
gives 
\begin{eqnarray*}
	\Ex \|S\|^2 
	& \le & v^2 \Biggl\{ 2\log d + 2\log(c/2) 
		- \log(2\log (dc/2)) \\
	&& \qquad \qquad + \ \
		\frac{8\sqrt{2\log(cd/2)} }{3\sqrt{2\log\big(\frac{cd}{2\sqrt{2\log(cd/2)}}\big)}  
			+ \sqrt{2\log\big(\frac{cd}{2\sqrt{2\log(cd/2)}}\big) + 8}}  \Biggr\} \\
	& =: & v^2 \{ 2 \log d + h_2(d)\}
\end{eqnarray*}
where it is easily checked that $h_2(d) \le 3$ for all $d \ge 1$. 
Moreover $h_2(d)$ is negative for $d >4.13795*10^{10}$. 
This completes the proof of the upper bound in (\ref{RademacherTypeTwoBestBounds}).

To prove the lower bound for $c_d$ in (\ref{RademacherTypeTwoBestBounds}), 
we use the lower bound of 
\cite{LedouxTalagrand:91}, Lemma 6.9, page 157 (which is, in this form, due to 
\cite{MR688642}).
This yields
\begin{eqnarray}
	c_d^2 \ge \frac{\lambda}{1+\lambda} t_o^2 + 
		\frac{1}{1+\lambda} d \int_{t_o}^{\infty} 4 t (1 - \Phi(t)) \, dt 
	\label{BasicVdVW-lowerBound}
\end{eqnarray}
for any $t_o > 0$, where $\lambda = 2 d (1 - \Phi (t_o))$. By using Komatsu's 
lower bound (\ref{ineq:KMSW}), we find that 
\begin{eqnarray*}
	\int_{t_o}^\infty  t (1 - \Phi(t)) \, dt 
	& \ge & \int_{t_o}^\infty \frac{2 t}{t + \sqrt{t^2+4}} \phi(t) \, dt \\
	& \ge & \frac{2 t_o}{t_o + \sqrt{t_o^2+4}} \int_{t_o}^\infty \phi(t) \, dt \\
	& = & \frac{2}{1 + \sqrt{1 + 4/t_o^2}} (1 - \Phi (t_o)) .
\end{eqnarray*}
Using this lower bound in (\ref{BasicVdVW-lowerBound}) yields
\begin{eqnarray}
	c_d^2 
	& \ge & \frac{\lambda}{1+\lambda} t_o^2 + 
		\frac{1}{1+\lambda} d \frac{8}{1 + \sqrt{1 + 4/t_o^2}} (1 - \Phi (t_o)) 
		\nonumber \\
	& = & \frac{2d(1-\Phi (t_o))}{1+ 2d(1-\Phi (t_o))}
		\left \{ t_o^2 + \frac{4}{1 + \sqrt{1 + 4/t_o^2}} \right \} 
		\nonumber \\
	& \ge & \frac{ \frac{4d}{t_o + \sqrt{t_o^2 + 4}} \phi (t_o)}
	             {1 + \frac{4d}{t_o + \sqrt{t_o^2 + 4}} \phi (t_o)}
		\left \{ t_o^2 + \frac{4}{1 + \sqrt{1 + 4/t_o^2}} \right \} .
	\label{LowerBoundAfterKomatsuTwice}
\end{eqnarray}
Now we let $c \equiv \sqrt{2 /\pi}$ and $\delta > 0$ and choose 
$$
	t_o^2 = 2 \log \left ( \frac{cd}{(2 \log (cd))^{(1+\delta)/2}} \right ).
$$
For this choice we see that $t_o \rightarrow \infty$ as $d \rightarrow \infty$,
\begin{eqnarray*}
	4 d \phi (t_o) = \frac{2d}{\sqrt{2\pi}} \cdot \frac{(2 \log (cd))^{(1+\delta)/2}}{cd} 
	= 2( 2 \log (cd))^{(1+\delta)/2} ,
\end{eqnarray*}
and 
\begin{eqnarray*}
	\frac{4 d \phi (t_o)}{t_o}
	= \frac{ 2( 2 \log (cd) )^{(1+\delta)/2} }{\{2 \log ( cd/ ( 2 \log (cd))^{(1+\delta)/2} )\}^{1/2} }
	\to \infty  
\end{eqnarray*}
as $d \to \infty$, so the first term on the RHS of (\ref{LowerBoundAfterKomatsuTwice}) 
converges to $1$ as $d \rightarrow \infty$, and it can be rewritten as 
\begin{eqnarray*}
	\lefteqn{ A_d \left \{t_o^2 + \frac{4}{1 + \sqrt{1 + 4/t_o^2}} \right \} } \\
	& = & A_d  \left \{ 2 \log \left ( \frac{cd}{(2 \log (cd))^{(1+\delta)/2}} \right )
		+  \frac{4}{1 + \sqrt{1 + 4/t_o^2}} \right \} \\
	& \sim & 1 \cdot \left \{ 2 \log d + 2 \log c - (1+\delta) \log (2 \log (cd)) +2 \right \} .
\end{eqnarray*}\\[-6ex]

To prove the upper bounds for $c_d$, we will use the 
upper bound of \cite{LedouxTalagrand:91}, Lemma 6.9, page 157 (which is, in this form, due to 
\cite{MR688642}).
For every $t_o>0$
\begin{eqnarray*}
c_{d}^2 \equiv \Ex \max_{1 \le j \le d} | Z_j|^2 
& \le & t_o^2 + d \int_{t_o}^\infty 2t P(|Z_1| > t) d t  \\
& = & t_o^2 + 4d \int_{t_o}^\infty t  (1 - \Phi (t)) dt \\
& \le & t_o^2 + 4d \int_{t_o}^\infty \phi (t) dt \qquad (\text{by Mills' ratio}) \\
& = & t_o^2 + 4d (1-\Phi (t_o)) .
\end{eqnarray*}
Evaluating this bound at $t_o = \sqrt{ 2 \log (d / \sqrt{2 \pi} ) } $ and then using Mills' ratio again yields
\begin{eqnarray}
c_{d}^2
& \le & 2 \log ( d / \sqrt{2 \pi} ) + 4 d (1 - \Phi (\sqrt{2 \log ( d / \sqrt{2 \pi} )} ) ) \nonumber  \\
& \le & 2 \log d - 2 \frac{1}{2} \log ( 2 \pi ) + 4 d \frac{\phi ( \sqrt{2 \log (d / \sqrt{2 \pi} )} ) }{\sqrt{2 \log ( d / \sqrt{2 \pi} )}} 
               \nonumber \\
& = & 2 \log d - \log ( 2\pi ) + \frac{2 \sqrt{2}}{\sqrt{ \log (d / \sqrt{2 \pi})}} 
            \label{AsymptoticFormWithFirstOrderCorrection}  \\
& \le & 2 \log d  \nonumber 
\end{eqnarray}
where the last inequality holds if
$$
      \frac{2\sqrt{2}}{\sqrt{ \log (d / \sqrt{2 \pi}) }} \le  \log (2\pi) , 
$$
or equivalently if
$$
      \log d \ \ge \ \frac{8}{(\log (2 \pi))^2} + \frac{\log(2\pi)}{2} \ = \ 3.28735  ...  ,
$$
and hence if $d \ge 27  > e^{3.28735...} \deq 26.77 $. The claimed inequality is easily verified 
numerically for $d =3, \ldots , 26$.  (It fails for $d=2$.)
As can be seen from (\ref{AsymptoticFormWithFirstOrderCorrection}), $2 \log d - \log (2 \pi)$ 
gives a reasonable approximation to $\Ex \max_{1 \le j \le d} Z_j^2$ for large $d$.
Using the upper bound in (\ref{ineq:KMSW}) instead of the second application of Mills' ratio and
choosing $t_o^2 = 2 \log (cd / \sqrt{2 \log (cd)} )$ 
with $c \ := \ \sqrt{2/\pi}$ yields the third bound for $c_d$ in (\ref{RademacherTypeTwoBestBounds})
with 
\begin{eqnarray*}
h_3 (d)
& = & - \log (\pi) - \log ( \log (cd))  \\
&& \ \  + \ \frac{8}{3 \sqrt{1- \frac{\log (2 \log (cd))}{2 \log (cd)}} + \sqrt{1- \frac{\log (2 \log (cd))}{2 \log (cd)} + \frac{4}{\log (cd)} }} .
\end{eqnarray*}
\strut\hfill	$\Box$

\subsection{Proofs for Section~\ref{sec: Bernstein}}

\paragraph{Proof of Lemma~\ref{lem: linexp}.}
It follows from $\Ex Z = 0$, the Taylor expansion of the exponential function 
and the inequality $\Ex |Z|^m \le \sigma^2 \kappa^{m-2}$ for $m \ge 2$ that
\begin{eqnarray*}
	\lefteqn{ \Ex \exp \Bigl( \frac{Z}{\kappa L} \Bigr)
		\ = \ 1 + \Ex \Bigl\{ \exp \Bigl( \frac{Z}{\kappa L} \Bigr) - 1 - \frac{Z}{\kappa L}
			\Bigr\} } \\
	& \le & 1 + \sum_{m=2}^{\infty} \frac{1}{m!} \frac{\Ex |Z|^m}{(\kappa L)^m}
		\ \le \ 1 + \frac{\sigma^2}{\kappa^2} \sum_{m=2}^{\infty} \frac{1}{m!} \frac{1}{L^m}
		\ = \ 1 + \frac{\sigma^2 {\rm e}(L)}{\kappa^2} .
\end{eqnarray*}\\[-5ex]
\strut\hfill$\Box$

\paragraph{Proof of Lemma~\ref{lem: Bernstein}.}
Applying Lemma~\ref{lem: linexp} to the $j$-th components 
$X_{i,j}$ of $X_i$ and $S_{n,j}$ of $S_n$ yields for all $L>0$,
$$
	\Ex \exp \Bigl( \frac{\pm S_{n,j}}{\kappa L} \Bigr)
	\ = \ \prod_{i=1}^n \Ex \exp \Bigl( \frac{\pm X_{i,j}}{\kappa L} \Bigr)
	\ \le \ \prod_{i=1}^n \exp \Bigl( \frac{\Var(X_{i,j}) {\rm e}(L)}{\kappa^2} \Bigr)
	\ \le \ \exp \Bigl( \frac{\Gamma {\rm e}(L)}{\kappa^2} \Bigr) .
$$
Hence
$$
	\Ex \cosh \Bigl( \frac{\|S_n\|_{\infty}}{\kappa L} \Bigr)
	\ = \ \Ex \max_{1 \le j \le d} \cosh \Bigl( \frac{S_{n,j}}{\kappa L} \Bigr)
	\ \le \ \sum_{j=1}^d
		\Ex \cosh \Bigl( \frac{S_{n,j}}{\kappa L} \Bigr)
	\ \le \ d \exp \Bigl( \frac{\Gamma {\rm e}(L)}{\kappa^2} \Bigr) .
$$
As in the proof of Lemma~\ref{lem: Veraar} we conclude that
\begin{eqnarray*}
	\Ex \|S_n\|_\infty^2
	& \le & (\kappa L)^2 \biggl( \log \Bigl(
		2 \Ex \cosh \Bigl( \frac{\|S_n\|_{\infty}}{\kappa L} \Bigr) \Bigr) \biggr)^2 \\
	& \le & (\kappa L)^2 \Bigl( \log(2d) + \frac{\Gamma {\rm e}(L)}{\kappa^2} \Bigr)^2 \\
	& = & \Bigl( \kappa L \log(2d) + \frac{\Gamma L \, {\rm e}(L)}{\kappa} \Bigr)^2 ,
\end{eqnarray*}
which is equivalent to the inequality stated in the lemma. \hfill $\Box$

\paragraph{Proof of Theorem~\ref{thm: Bernstein-Nemirovski}.}
For fixed $\kappa_o > 0$ we split $S_n$ into $A_n + B_n$ as described before. 
Let us bound the sum $B_n$ first: For this term we have 
\begin{eqnarray*}
	\|B_n\|_{\infty}
	& \le & \sum_{i=1}^n \bigl\{ 1_{[\|X_i \|_\infty > \kappa_o ]} \|X_i\|_{\infty}
		+ \Ex(1_{[\|X_i \|_\infty > \kappa_o ]} \|X_i\|_\infty) \bigr\} \\
	& = & \sum_{i=1}^n \bigl\{ 1_{[\|X_i \|_{\infty} > \kappa_o ]} \|X_i\|_\infty
		- \Ex (1_{[\|X_i \|_{\infty} > \kappa_o ]} \|X_i\|_\infty) \bigr\} \\
	&& \qquad + \ 2 \sum_{i=1}^n \Ex ( 1_{[\|X_i \|_\infty > \kappa_o ]} \| X_i \|_{\infty} ) \\
	& =: & B_{n1} + B_{n2} .
\end{eqnarray*}
Therefore, since $\Ex B_{n1} = 0$,
\begin{eqnarray*}
	\Ex \|B_n\|_\infty^2 
	& \le & \Ex (B_{n1} + B_{n2})^2 \ = \ \Ex B_{n1}^2 + B_{n2}^2 \\
	& = & \sum_{i=1}^n \Var \bigl( 1_{[\|X_i \|_\infty > \kappa_o]} \|X_i\|_\infty \bigr)
		+ 4 \Bigl( \sum_{i=1}^n
			\Ex ( \| X_i \|_{\infty} 1_{[\|X_i \|_{\infty} > \kappa_o ]} ) 	\Bigr)^2 \\
	& \le & \sum_{i=1}^n \Ex \|X_i\|_\infty^2
		+ 4 \Bigl( \sum_{i=1}^n  \frac{ \Ex \| X_i \|_{\infty}^2 }{\kappa_o} \Bigr)^2 \\
	& = & \Gamma + 4 \frac{\Gamma^2}{ \kappa_o^2} ,
\end{eqnarray*}
where we define $\Gamma := \sum_{i=1}^n \Ex \|X_i\|_\infty^2$.

The first sum, $A_n$, may be bounded by means of Lemma~\ref{lem: Bernstein} 
with $\kappa = 2\kappa_o$, utilizing the bound
$$
	\Var(X_{i,j}^{(a)})
	\ = \ \Var \bigl( 1_{[\|X_i\|_\infty \le \kappa_o]} X_{i,j} \bigr)
	\ \le \ \Ex \bigl( 1_{[\|X_i\|_\infty \le \kappa_o]} X_{i,j}^2 \bigr)
	\ \le \ \Ex \|X_i\|_\infty^2 .
$$
Thus
$$
	\Ex \|A_n\|_\infty^2
	\ \le \ \Bigl( 2\kappa_o L \log (2d) + \frac{\Gamma L \, {\rm e}(L)}{2\kappa_o} \Bigr)^2 .
$$

Combining the bounds we find that 
\begin{eqnarray*}
	\sqrt{\Ex \| S_n \|_\infty^2} 
	& \le & \sqrt{\Ex \|A_n\|_\infty^2} + \sqrt{\Ex \|B_n\|_\infty^2 } \\
	& \le & 2\kappa_o L \log (2d) + \frac{\Gamma L {\rm e}(L)}{2\kappa_o} 
		+ \ \sqrt{\Gamma} + 2 \frac{\Gamma}{\kappa_o} \\
	& = & \alpha \kappa_o + \frac{\beta}{\kappa_o} + \sqrt{\Gamma} ,
\end{eqnarray*} 
where $\alpha := 2L \log(2d)$ and $\beta := \Gamma (L \, {\rm e}(L) + 4) / 2$.
This bound is minimized if $\kappa_o = \sqrt{\beta/\alpha}$ with minimum value
$$
	2 \sqrt{\alpha\beta} + \sqrt{\Gamma}
	\ = \ \bigl( 1 + 2 \sqrt{L^2 {\rm e}(L) + 4L} \sqrt{\log(2d)} \bigr) \sqrt{\Gamma} ,
$$
and for $L = 0.407$ the latter bound is not greater than
$$
	\bigl( 1 + 3.46 \sqrt{\log(2d)} \bigr) \sqrt{\Gamma} .
$$

In the special case of symmetrically distributed random vectors $X_i$, 
our treatment of the sum $B_n$ does not change, but in the bound for $\Ex \|A_n\|_\infty^2$ 
one may replace $2\kappa_o$ with $\kappa_o$, because $\Ex X_i^{(a)} = 0$. Thus
\begin{eqnarray*}
	\sqrt{\Ex \| S_n \|_\infty^2} 
	& \le & \kappa_o L \log (2d) + \frac{\Gamma L {\rm e}(L)}{\kappa_o} 
		+ \ \sqrt{\Gamma} + 2 \frac{\Gamma}{\kappa_o} \\
	& = & \alpha' \kappa_o + \frac{\beta'}{\kappa_o} + \sqrt{\Gamma}
		\qquad	\bigl( \text{with} \
			\alpha' := L \log(2d), \beta' := \Gamma (L \, {\rm e}(L) + 2) \bigr) \\
	& = & \bigl( 1 + 2 \sqrt{L^2 {\rm e}(L) + 2L} \sqrt{\log(2d)} \bigr) \sqrt{\Gamma}
		\qquad	\bigl(\text{if} \ \kappa_o = \sqrt{\beta'/\alpha'} \bigr) .
\end{eqnarray*} 
For $L = 0.5$ the latter bound is not greater than
$$
	\bigl( 1 + 2.9 \sqrt{\log(2d)} \bigr) \sqrt{\Gamma} .
	\eqno{\Box}
$$

\paragraph{Acknowledgements.}
The authors owe thanks to the referees for a number of suggestions which 
resulted in a considerable improvement in the article.
The authors are also grateful to Ilya Molchanov for drawing their attention to 
Banach-Mazur distances, and to Stanislaw Kwapien and Vladimir Koltchinskii 
for pointers concerning type and co-type proofs and constants.
This research was initiated during the opening week of the program on 
``Statistical Theory and Methods for Complex,
High-Dimensional Data'' held at the {\sl Isaac Newton Institute for Mathematical Sciences} 
from 7 January to 27 June, 2008, and was made possible in 
part by the support of the Isaac Newton Institute for visits of various periods by 
D\"umbgen, van de Geer, and Wellner.
The research of Wellner was also supported in part by NSF grants DMS-0503822 
and DMS-0804587. The research of D\"{u}mbgen and van de Geer was 
supported in part by the Swiss National Science Foundation.


\bibliographystyle{monthly}
\bibliography{Nemirovski-new.bib}

\bigskip

\noindent\textbf{Lutz D\"{u}mbgen} received his Ph.D. from 
Heidelberg University in 1990. From 1990-1992 he was a Miller 
research fellow at the University of California at Berkeley. Thereafter he 
worked at the universities of Bielefeld, Heidelberg and L\"{u}beck. Since 2002 he is 
professor of statistics at the University of Bern. 
His research interests are nonparametric, multivariate and computational statistics.

\noindent\textit{Institute of Mathematical Statistics and Actuarial Science, 
University of Bern, Alpenegg\-strasse 22, CH-3012 Bern, Switzerland\\
duembgen@stat.unibe.ch}

\bigskip

\noindent\textbf{Sara A. van de Geer} obtained her Ph.D. at Leiden 
University in 1987. She worked at the Center for Mathematics and 
Computer Science in Amsterdam, at the Universities of Bristol, Utrecht, Leiden 
and Toulouse, and at the Eidgen\"{o}ssische Technische Hochschule in 
Z\"{u}rich (2005-present). Her research areas are empirical processes, 
statistical learning, and statistical theory for high-dimensional data. 

\noindent\textit{Seminar for Statistics, ETH Zurich, CH-8092 Zurich, Switzerland\\
geer@stat.math.ethz.ch}

\bigskip

\noindent\textbf{Mark C. Veraar}  received his Ph.D. from Delft University 
of Technology in 2006. In the year 2007 he stayed as a PostDoc 
researcher in the European RTN project ``Phenomena in High Dimensions'' 
at the IMPAN institute in Warsaw (Poland). In 2008 he spent one year as an 
Alexander von Humboldt fellow at the University of Karlsruhe (Germany). 
Since 2009 he is Assistant Professor at Delft University of Technology (the Netherlands). 
His main research areas are probability theory, partial differential equations and functional analysis.

\noindent\textit{Delft Institute of Applied Mathematics,
Delft University of Technology, P.O. Box 5031, 2600 GA Delft, The Netherlands\\
m.c.veraar@tudelft.nl, mark@profsonline.nl} 

\bigskip

\noindent\textbf{Jon A. Wellner} received his B.S. from the University of Idaho in 1968 
and his Ph.D. from the 
University of Washington in 1975.  
He got hooked on research in 
probability and statistics during graduate school at the UW in the early 1970's, 
and has enjoyed both teaching
and research at the University of Rochester (1975-1983) and the 
University of Washington (1983-present).  
If not for probability theory and statistics, he might be a ski bum.

\noindent\textit{Department of Statistics, Box 354322, University of Washington, 
Seattle, WA  98195-4322\\
jaw@stat.washington.edu}

\end{document}